\documentclass[12pt, oneside]{article}   	
\usepackage{amsmath}
\usepackage{geometry} 
\usepackage{mathtools}
\usepackage{amsthm}  
\usepackage{subfig}

\newcommand{\bM}{{\bf{M}}}
\newcommand{\bH}{{\bf{H}}}

\newcommand{\br}{{\bf{x}}}
\newcommand{\by}{{\bf{y}}}
\newcommand{\bn}{{\bf{n}}}
\newcommand{\be}{{\bf{e}}}
\newcommand{{\bome}}{{\boldsymbol{\omega}}}
\newcommand{\te}{{\tilde{e}}}
\newcommand{\brp}{{\bf{x}^{\prime}}}
\newcommand{\bnp}{{\bf{n}^{\prime}}}
\newcommand{\bx}{{\bf{x}}}
\newcommand{\bz}{{\bf{z}}}
\newcommand{\bxp}{{\bf{x}^{\prime}}}

\newtheorem{Remark}{Remark}
\newtheorem{Theorem}{Theorem}
\newtheorem{Lemma}{Lemma}
\newtheorem*{Proof}{Proof}
\usepackage[nottoc]{tocbibind}
\usepackage{authblk}
\usepackage{stmaryrd}

\usepackage{graphicx}				
\usepackage{amssymb}

\title{A hybrid boundary integral-PDE approach for the approximation of the demagnetization potential in micromagnetics}

\author[1]{Doghonay Arjmand}
\author[2]{V\'ictor Mart\'inez Calzada}
\affil[1]{Department of Information Technology, Division of Scientific Computing, Uppsala University}
\affil[2]{Departamento de Ingenier\'ia Aeron\'autica, Universidad Polit\'ecnica Metropolitana de Hidalgo, Hidalgo, M\'exico}

\date{}							

\begin{document}
\maketitle

\begin{abstract}
The demagnetization field in micromagnetism is given as the gradient of a potential which solves a partial differential equation (PDE) posed in $\mathbb{R}^d$. In it's most general form, this PDE is supplied with continuity condition on the boundary of the magnetic domain and the equation includes a discontinuity in the gradient of the potential over the boundary. Typical numerical algorithms to solve this problem relies on the representation of the potential via the Green's function, where a volume and a boundary integral terms need to be accurately approximated. From a computational point of view, the volume integral dominates the computational cost and can be difficult to approximate due to the singularities of the Green's function. In this article, we propose a hybrid model, where the overall potential can be approximated by solving two uncoupled PDEs posed in bounded domains, whereby the boundary conditions of one of the PDEs is obtained by a low cost boundary integral. Moreover, we provide a convergence analysis of the method under two separate theoretical settings; periodic magnetisation, and high-frequency magnetisation. Numerical examples are given to verify the convergence rates. 
\end{abstract} 

\section{Introduction}
Magnetic behaviour of materials can be modelled at different time and length scales. At the smallest scale, the spin and orbital movements of magnetic moments can be modelled by electronic structure calculations. Despite the accuracy that the electronic calculations can provide, the resolution of full electronic structure is prohibitively expensive even with the best available computers today. At a slightly larger (atomistic) scale, the magnetic behaviour can be studied by modelling the interaction of atomic moments, whereby certain parameters (e.g., exchange coefficient between atoms) maybe upscaled using local electronic calculations. The magnetic behaviour at the atomistic scale is described by a set of coupled nonlinear ordinary differential equations known as Landau-Lifshitz-Gilbert equations, and the atomic moments interact with each other due to various short and long range interactions as well as external excitations via an applied magnetic field. Micromagnetism, in particular, deals with the study of magnetic materials below micrometer temporal and spatial scales. In this regime, magnetization is governed by the following (still named as Landau-Lifshitz-Gilbert) nonlinear partial differential equation (PDE):
\begin{equation}\label{eqn_LLG}
    \partial_t \bM(t,\bx) = - \bM \times \bH  - \gamma \bM \times \left( \bM \times \bH \right), \quad |\bM| = 1, 
\end{equation}
where
\begin{equation*}
    \bH(t,\bx,\bM) = \bH_{ex}(\bM) + \bH_{ext}(t,\bx) + \bH_{ani}(\bx) + \bH_{str}(\bM) + \bH_{dem}(\bM). 
\end{equation*}
Here $\bH_{ex}(\bM),\bH_{ext}(t,\bx),\bH_{ani}(\bx),\bH_{str}(\bM),$ and $\bH_{dem}(\bM)$ represent the exchange, external, anisptropy, magnetorestrictive, and the demagnetization fields respectively, see e.g., \cite{Cervera2007,Cimrak2008} for a detailed description and numerical treatment of these terms. Important for the aim of this paper is the demagnetization field $\bH_{dem}$, which is originally due to long range interactions of atomic moments. Computing $\bH_{dem}$ at an atomistic scale is prohibitively expensive for problem sizes seen in engineering applications, and a continuum approximation is often preferred. An efficient computation of $\bH_{dem}$ is not only important for solving the continuum model \eqref{eqn_LLG} but also extremely vital in maintaining a low computational cost in multiscale algorithms, such as \cite{ARJMAND201999,Arjmand2017AtomisticcontinuumMM,Arjmand2020LongRange}, where $\bH_{dem}$ is computed on a continuum level and used in the atomistic models to decrease the overall cost of the simulations; see also \cite{GarciaSanchez2005,Jourdan2008,Kakay_etal_2014,Poluektov2018,DeLucia2016,Leitenmaier_Runborg_2022} for other multiscale strategies in micromagnetism.

\subsection*{Continuum description of $\bH_{dem}$:}

Evident from equation \eqref{eqn_LLG}, simulations of magnetic materials at the micromagnetic regime require the solution of the demagnetization field $\bH_{dem}(\bM) \in \mathbb{R}^3$, for all $t \in \mathbb{R}^{+},\bx \in \mathbb{R}^{d},$ and $\bM(t,\bx) \in \mathbb{R}^{d} $, and due to it's long range nature this is the most time consuming part in micromagnetic simulations. The field $\bH_{dem}$ is given as $\bH_{dem} = -\nabla u$, where 

\begin{align} \label{Eqn_Demag}
&-\triangle u(\bx)  = 
\begin{cases} - \nabla \cdot \bM(\br) & \text{ if } \br \in \Omega  \\
0 & \text{ if } \br \in \mathbb{R}^3 / \bar{\Omega},
\end{cases}
\\
&u \text{ is continuous on } \partial\Omega \nonumber \\
&\llbracket \bn \cdot \nabla u \rrbracket = \bM \cdot \bn \text{ on } \partial{\Omega}.  \nonumber
\end{align}
Here $\Omega \subset \mathbb{R}^3$ is an open bounded domain filled with a magnetic material, $\bn$ is the unit normal vector taken to be positive in the outward direction to the boundary $\partial \Omega$ of the domain $\Omega$, and $\llbracket \cdot \rrbracket$ denotes a jump accross the interface between the magnetic material and the rest of $\mathbb{R}^d$. It is also assumed that the potential $u$ along with it's gradient decays to zero as $\br \to \infty$. In simulations the infinite domain problem \eqref{Eqn_Demag} needs to be solved many times as the magnetization vector $\bM$ changes over time\footnote{Here we omit the temporal dependency of $\bM$ as it does not serve the purpose of this work}. Naive approximations based on finite element method (FEM) or finite differences (FD) will suffer from the influence of the artificial boundary conditions, due to the need to truncate the infinite domain. Moreover, these approaches result in very large linear systems in three dimensions, which is prohibitive from a computational viewpoint. 

One popular solution method is to represent the solution $u$ via the following integral representation
\begin{equation*}
u(\br) = - \int_{\Omega} G(\br,\brp) \nabla^{\prime} \cdot \bM(\brp) \; d\br^{\prime}  + \int_{\partial\Omega} G(\br,\brp) \bM(\brp) \cdot \bnp \; d\brp, 
\end{equation*}
where $G$ denotes the infinite domain Green's function for the Laplace operator. A derivation of this integral representation is provided in Section \ref{Sec_Formal_Integral}. Moreover, equivalent formulations can be derived by integrating by part and moving the derivative to the Green's function. However, this is not favourable due to strong singularities of the gradient of the Green's function. From a computational viewpoint, the main issue with this integral representation is that the cost of computing the solution inside the domain $\Omega$ scales as $O(N^{2d} + N^{2d-1})$, where $N$ denotes the number of degrees of freedom in one direction. This is computationally very expensive, in particular for three dimensional problems and more efficient approaches are needed.

In this article, we develop and analyse a method for reducing the computational cost to $O(N^d + N^{2d-2})$. The proposed approach here relies on splitting the solution into two simple uncoupled (and therefore easily parallelizable) PDE problems over bounded domains, whereby the boundary conditions of one of the PDEs is given by a boundary integral formula. There are conceptually similar approaches in the literature, \cite{Fredkin_Koehler_1990,Cervera_Roma_2006}, which we discuss in the following subsection.

\subsection*{Existing methods for computing $u$}
\label{Sec_Existing_Methods}
In the literature, there are essentially two sets of efficient methods to compute $u$. The first approach is to design fast multipole algorithms in order to accelerate the computations of the integral formulations, see e.g., \cite{Long_2006,PALMESI2017409,EXL2014490} and the references therein. Another class of algorithms are based on hybrid PDE-boundary integral methods (BEM), see \cite{Fredkin_Koehler_1990,Cervera_Roma_2006}, which is conceptually similar to the idea of the present paper. The main difference between the approaches of Fredkin-Koehler \cite{Fredkin_Koehler_1990} and Cervera-Roma \cite{Cervera_Roma_2006} is the fact that the former uses a double layer potential formulation for the boundary integral, which means that the boundary integral includes the gradient of the Green's function,  while the latter uses a single layer potential potential formulation, which includes the Green's function itself in the boundary integral. Thus the latter is more favourable from a computational point of view since the strong singularities in the gradient of the Green's function are not present in the algorithm. Now, we present the hybrid algorithm from  \cite{Cervera_Roma_2006} as it bears similarities with the methodology proposed in this paper. The idea of \cite{Cervera_Roma_2006} is to split the solution $u = v + w$, so that 
\begin{align}
    \label{Eqn_Cervera_Roma_w}
   - \Delta w(\bx) &= -\nabla \cdot \bM(\bx), \quad \bx \in \Omega, \\
            w(\bx) &= 0, \quad \text{ on } \quad \partial\Omega, \nonumber
\end{align}
where $w$ is extended to zero in $\mathbb{R}^d \backslash \Omega$, and  
\begin{align}
    \label{Eqn_Cervera_Roma_v}
   - \Delta v(\bx) &= 0, \quad \bx \in \Omega \cup \{ \mathbb{R}^d \backslash \overline{\Omega} \}, \\
            \llbracket v \rrbracket  &= 0, \quad \text{ on } \quad \partial\Omega, \nonumber \\
            \llbracket \bn \cdot \nabla v \rrbracket  &= -\bM \cdot \bn  + \bn \cdot \nabla w, \quad \text{ on } \quad \partial\Omega. \nonumber 
\end{align}
The authors then express the solution $v$ in terms of the following boundary integral 
\begin{align}    \label{Eqn_Cervera_Roma_Boundary}
    v(\bx) = \int_{\partial \Omega} G(\bx,\by) g(\by) \; d\by,
\end{align}
where $g(\by) = -\bM(\by) \cdot \bn  + \bn \cdot \nabla w(\by)$. Note that the cost of computing the solution $v$ via \eqref{Eqn_Cervera_Roma_Boundary} is $O(N^{2d-1})$. To decrease the cost of the algorithm, the authors use \eqref{Eqn_Cervera_Roma_Boundary} to compute $v$ only on the boundary $\partial \Omega$, and then solve the Laplace equation \eqref{Eqn_Cervera_Roma_v} with Dirichlet boundary conditions coming from the boundary integral \eqref{Eqn_Cervera_Roma_Boundary}. Assuming a linearly scaling PDE solver, e.g., a multigrid method, the cost is then reduced to $O(N^d + N^{2d-2})$. The method proposed in this paper has the same cost as in the above-mentioned algorithm, and it also relies on a single layer potential formulation for the boundary integral part, see Section \ref{Sec_HybridApproach}. The main difference is that in the modelling strategy above, the PDE solutions $v$ and $w$ are coupled through the boundary conditions of $v$, whereas the splitting used in this paper results in two totally uncoupled problems; which is desirable for parallel computations (considering the parallel computing capabilities of computers today).

This article is structured as follows: In Section \ref{Sec_Prel}, we introduce the notations and preliminaries. In Section \ref{Sec_Formal_Integral}, we present a derivation of a well-known integral representation for $u$. The formulation of the proposed hybrid approach is presented in Section \ref{Sec_HybridApproach}, and it's analysis is included in Section \ref{Sec_Analysis}. Finally, in Section \ref{Sec_Num_Res}, we provide numerical examples in two dimensions to verify our theoretical findings.

\section{Preliminaries}
\label{Sec_Prel}
\begin{itemize}

    \item[{\bf 1.}] Throughout the sequel, we denote the vector valued quantities by bold face letters; e.g., $\bx$ represents a vector in $\mathbb{R}^{d}$, while $u$ would represent a scalar valued function. 

    \item[{\bf 2.}] Let $K = [0,1]^d$, we denote the space of periodic functions $f$ belonging to $W^1_{per}(K)$ as 
    \begin{align*}
        W^1_{per}(K):= \{ f \text{ is } K\text{-periodic}: \int_{K} f \; dx  = 0,  \text{ and } \| f \|_{L^2(K)}  + \|  \nabla f \|_{L^2(K)} \leq C \},
    \end{align*}
    where $C$ will denote a generic constant whose dependence on critical parameters are explicitly stated in the manuscript. 
    \item[{\bf 3.}] We denote the class $\mathcal{H}_{\omega_0}$ of $L^2(\mathbb{R}^d)$ integrable functions with vanishing Fourier modes in the origin as 
\begin{align*}
    \mathcal{H}_{\omega_0} := \{ f \in L^2(\mathbb{R}^d): \hat{f}(\bome) = 0, \quad \forall\;|\bome|_2 \leq \omega_0, \quad \text{ for some } \omega_0 >0 \}
\end{align*}

    \item[{\bf 4.}] Let $K = [0,1]^d$, and consider the Laplace operator $-\Delta$ equipped with periodic boundary conditions on $K$. Then the eigenvalues $\{ \lambda_j \}_{j=0}^{\infty}$ of $-\Delta$ are non-negative, see \cite{Safarov_Vassiliev_1997}, and assumed to be ordered as follows 
    \begin{align*}
        0 = \lambda_0 < \lambda_1 \leq \lambda_2 \leq \lambda_3 \leq \ldots
    \end{align*}
    Important for the upcoming estimates is the smallest positive eigenvalue $\lambda_1 \geq \frac{\pi^2}{d}$. This is due to the fact that $\lambda_1  = C_p(K)^{-2}$, where $C_p$ is the Poincaré-Wirtinger constant for the domain $K$ and is bounded by $C_p \leq \frac{diam(K)}{\pi} = \frac{\sqrt{d}}{\pi}$, see e.g., \cite{Payne_Weinerberger_1960}. Moreover, the corresponding eigenfunctions $\{ \varphi_j \}_{j=0}^{\infty}$ are orthonormal in $L^2(K)$ and satisfy the following properties
    \begin{itemize}
\item[P1.] $L^2(K)$ orthogonality:
\begin{align*}
        \langle  \varphi_i,\varphi_j \rangle_{K}:= \int_K \varphi_i(\bx) \varphi_j(\bx) \; d\bx = \begin{cases} 0 & \quad \text{if} \quad i\neq j \\
        1 & \quad \text{if} \quad i=j.
        \end{cases}
    \end{align*}
\item[P2.] Zero average property:
   \begin{align*}
        \langle  \varphi_i, 1 \rangle_{K} = \begin{cases} 1 & \quad \text{if} \quad i= 0, \\
        0 & \quad \text{if} \quad i\geq 1.
        \end{cases}
    \end{align*}
\end{itemize}

\end{itemize}

\section{Derivation of a known integral formulation}
\label{Sec_Formal_Integral}
In the literature, see e.g., \cite{aharoni2000introduction}, an integral representation of the potential $u$ is given as

\begin{equation} \label{Eqn_Integral_Representation}
u(\br) = - \int_{\Omega} G(\br,\brp) \nabla^{\prime} \cdot \bM(\brp) \; d\br^{\prime}  + \int_{\partial\Omega} G(\br,\brp) \bM(\brp) \cdot \bnp \; d\brp, 
\end{equation}
where the Green's function $G(\br,\brp)$ satisfies
\begin{equation} \label{Eqn_Greens_Function}
-\triangle G(\br,\brp) = \delta(\br,\brp), \text{ in } \mathbb{R}^d.
\end{equation}
It has come to the attention of the authors that, in the literature, sometimes the second integral term in \eqref{Eqn_Integral_Representation} is ignored. However, a complete description of the potential $u$ needs to accommodate both terms. Since we do not find a mathematical derivation of this formula in the literature, and for the sake of completeness, we provide a derivation of this integral representation. Note that other equivalent formulations are possible to derive by further integrating by part the first integral in the right hand side and moving the derivative to the Green's function but this is often computationally unfavourable due to the strong singularities that the gradient possesses.   
\begin{Lemma} Suppose $u$ satisfies equation \eqref{Eqn_Demag} along with the boundary conditions, then $u$ satisfies the integral representation \eqref{Eqn_Integral_Representation}.
\end{Lemma}
\begin{Proof} Let 
\begin{equation*}
-\Delta u(\brp)  = F(\brp), \text { where }  F(\brp) = \begin{cases} -\nabla \cdot \bM(\brp), & \brp \in \Omega_1 \\
0, & \brp \in \Omega_2,
\end{cases}
\end{equation*}
where $\Omega_1 = \Omega$ and $\Omega_2  = \mathbb{R}^{d} / \bar{\Omega}_1$. Then it follows that 
\begin{align} \label{Eqn_Aux}
\int_{\mathbb{R}^{d}} -\Delta u(\brp) G(\br,\brp) \; d\brp = \int_{\mathbb{R}^{d}} F(\brp) G(\br,\brp) \; d\brp = -\int_{\Omega_1} G(\br,\brp) \nabla \cdot \bM(\brp) \; d\brp.
\end{align}
 Moreover, we denote by $\bn_1$ and $\bn_2$ the unit normal to $\Omega_1$ and $\Omega_2$ chosen in the outward direction to respective domains so that $\bn_1=-\bn_2$, and we now focus only on the left hand side.
\begin{align*}
\int_{\mathbb{R}^{d}} -\Delta u(\brp) G(\br,\brp) \; d\brp & = \int_{\Omega_1} -\Delta u(\brp) G(\br,\brp) \; d\brp + \int_{\Omega_2} -\Delta u(\brp) G(\br,\brp) \; d\brp \\
&= \int_{\Omega_1} \nabla u(\brp) \cdot \nabla G(\br,\brp) \; d\brp - \int_{\partial\Omega_1} \nabla u(\brp) \cdot \bnp_1 G(\br,\brp) \; d\brp \\
&+\int_{\Omega_2} \nabla u(\brp) \cdot \nabla G(\br,\brp) \; d\brp - \int_{\partial\Omega_2} \nabla u(\brp) \cdot \bnp_2 G(\br,\brp) \; d\brp \\
&= -\int_{\Omega_1} u(\brp)  \Delta G(\br,\brp) \; d\brp +  \int_{\partial\Omega_1}  u(\brp)  \nabla G(\br,\brp) \cdot \bnp_1  \; d\brp \\
& -\int_{\Omega_2} u(\brp) \Delta G(\br,\brp) \; d\brp +  \int_{\partial\Omega_2} u(\brp)  \nabla G(\br,\brp) \cdot \bnp_2  \; d\brp \\
& - \int_{\partial\Omega_1} \nabla u(\brp) \cdot \bnp_1 G(\br,\brp) \; d\brp - \int_{\partial\Omega_2} \nabla u(\brp) \cdot \bnp_2 G(\br,\brp) \; d\brp .
\end{align*}
Note that in the above, we also used the fact that the solution $u(\br)$ along with it's gradient decays to zero as $\br \to \infty$. Now using $\partial\Omega_1  = \partial\Omega_2$, and $\bnp_1 = -\bnp_2$, the continuity of the solution $u$ on $\partial\Omega$, as well as equation \eqref{Eqn_Greens_Function}, and the jump condition $\llbracket \bn \cdot \nabla u \rrbracket = \bf{M} \cdot \bf{n}$  on the boundary $\partial \Omega$ we obtain
\begin{align*}
\int_{\mathbb{R}^{d}} -\Delta u(\brp) G(\br,\brp) \; d\brp & =  -\int_{\Omega_1} u(\brp) \Delta G(\br,\brp) \; d\brp -\int_{\Omega_2} u(\brp) \Delta G(\br,\brp) \; d\brp \\
&- \int_{\partial\Omega_1}\llbracket \bn \cdot \nabla u \rrbracket G(\br,\brp) \; d\brp \\
&= \int_{\Omega_1} u(\brp) \delta(\br,\brp) \; d\brp +\int_{\Omega_2} u(\brp) \delta(\br,\brp) \; d\brp -  \int_{\partial\Omega_1}\bM \cdot \bn_1 G(\br,\brp) \; d\brp
\end{align*}
On the other hand, \eqref{Eqn_Aux} together with the above formula yields
\begin{align*}
 \int_{\Omega_1} u(\brp) \delta(\br,\brp) \; d\brp +\int_{\Omega_2} u(\brp) \delta(\br,\brp) \; d\brp -  \int_{\partial\Omega_1}\bM \cdot \bn_1 G(\br,\brp) \; d\brp = - \int_{\Omega_1} G(\br,\brp) \nabla \cdot \bM(\brp) \; d\brp
\end{align*}
If $\br \in \Omega_1$, then the second term in the left hand side vanishes and we obtain
\begin{align*}
u(\br) =  \int_{\partial\Omega_1}\bM(\brp) \cdot \bn_1 G(\br,\brp) \; d\brp -  \int_{\Omega_1} G(\br,\brp) \nabla \cdot \bM(\brp) \; d\brp
\end{align*}
Similarly if $\br \in \Omega_2$, then the first term in the left hand side vanishes and we get
\begin{align*}
u(\br) =  \int_{\partial\Omega_1}\bM(\brp) \cdot \bn_1 G(\br,\brp) \; d\brp -  \int_{\Omega_1} G(\br,\brp) \nabla \cdot \bM(\brp) \; d\brp
\end{align*}
Therefore, for any $\br \in \mathbb{R}^{d}$, we have 
\begin{align} \label{Eqn_Formal_Formulation}
u(\br) =  \int_{\partial\Omega_1}\bM(\brp) \cdot \bn_1 G(\br,\brp) \; d\brp -  \int_{\Omega_1} G(\br,\brp) \nabla \cdot \bM(\brp) \; d\brp.
\end{align}
\end{Proof}


\section{A Hybrid Approach}
\label{Sec_HybridApproach}
The classical integral representation \eqref{Eqn_Integral_Representation} of the solution $u$ includes a volume and a boundary integral term. The cost of this integral computation is dominated by the volume integral, while the boundary integral, which is as a result of the jump condition in the normal derivative of the solution, can be computed more efficiently. The idea behind our hybrid approach is to approximate both the volume integral and the boundary integral (for interior points) using two uncoupled elliptic PDEs, whereby the boundary conditions of the second PDE is obtained by a boundary integral. In particular, we approximate the representation \eqref{Eqn_Integral_Representation} by the following splitting
\begin{equation} \label{Eqn_Hybrid_Representation}
u(\bx) \approx v_{T,R}(\bx)  + b(\bx), 
\end{equation}
where
\begin{align} \label{Eqn_RegularisedElliptic}
-\Delta v_{T,R}(\bx)&= F - e^{T  \Delta } F, \quad \bx \in K_R, \\
v_{T,R}(\bx)&= 0, ´\quad \text{on} \quad \partial K_R,  \nonumber
\end{align}
and
\begin{align} \label{Eqn_b}
    -\Delta b(\bx) &= 0, \quad \bx \text { in } \Omega \\
    b(\bx ) &= \int_{\partial\Omega} G(\bx,\by) \bM(\by) \cdot \bn \; d\by, \quad \text{ on } \partial\Omega \nonumber
\end{align}
Here $R>1$ so that $\Omega \subset K_{R}$, and $T>0$ is a fixed parameter whose optimal value will be evident in upcoming analysis. If the term $v_{T,R}$ approximates the volume integral:
\begin{equation} \label{Eqn_Approximation}
 v_{T,R}(\bx) \approx - \int_{\Omega} G(\bx,\by) \nabla_{\by} \cdot \bM(\by) \; d\by =:w_{\infty}(\bx), 
\end{equation}
where $G$ is the Green's function for the whole space defined in  \eqref{Eqn_Greens_Function}, then the solution $u$ in \eqref{Eqn_Hybrid_Representation} and  \eqref{Eqn_Integral_Representation} will stay close to each other. This is due to the fact that, by construction, the term $b$ is exactly equal to the boundary integral term in \eqref{Eqn_Integral_Representation}. To see this assume that $\bx \in \Omega$, then we recover equation \eqref{Eqn_b} by simply considering
\begin{equation*}
    -\Delta_{\bx} \int_{\partial \Omega} G(\bx,\by) \bM(\by) \cdot \bn \; d\by = 0.
\end{equation*}
Therefore, interesting from an analysis point of view is the convergence of the term $v_{T,R}$ solving \eqref{Eqn_RegularisedElliptic} to the right hand side $w_{\infty}$ in \eqref{Eqn_Approximation} which solves the following PDE\footnote{The choice for the notation $w_{\infty}$ will be clear in a nutshell}
\begin{align} \label{Eqn_InfDomain_Elliptic}
-\Delta w_{\infty}(\bx) = F(\bx), \quad \bx \in \mathbb{R}^d,
\end{align}
where 
\begin{align*} 
F(\bx)  = \begin{cases} -\nabla  \cdot M(\bx) &  \text{ if }  \bx  \in  \Omega, \\
0 & \text{ if }  \bx  \in  \mathbb{R}^d / \Omega.
\end{cases}
\end{align*}
A detailed numerical study of the convergence $v_{T,R} \to w_{\infty}$ is given in the numerical results Section. 

One naive approach to approximate the solution to \eqref{Eqn_InfDomain_Elliptic} is to truncate the infinite domain to a finite domain and solve the bounded domain problem (equipped with homogeneous Dirichlet boundary conditions) using a standard discretization scheme. However, decay of the Green's function associated with problem \eqref{Eqn_InfDomain_Elliptic} is poor ($1/|\bx|^{d-2}$), and the artificial boundary conditions posed on the boundary of the truncated domain will significantly pollute the interior solution unless a very large computational geometry is used, see Theorem \ref{Thm_Analysis_1}. At this point, we resort to a recently developed method in the context of homogenization problems, see \cite{Abdulle_Arjmand_Paganoni_2023}, where the authors relate the solutions of the heat equation to that of elliptic equations. This approach is then much more efficient due to the Gaussian decay of the Green's function associated with the heat equation. This fact will allow for comping up with a PDE model posed over a relatively small domain (but still larger than the original geometry of interest where the magnetic material is located). To motivate the idea, let $w(t,\bx)$ solve the following heat equation in $\mathbb{R}^d$
\begin{align} \label{Eqn_InfDomain_Heat}
\partial_t w(t,\bx)-\Delta w(t,\bx)&= 0, \quad \bx \in \mathbb{R}^d, \\
w(0,\bx)  &= F(\bx). \nonumber
\end{align}
Then using the decaying properties of the solutions to the heat equation as $t \to \infty$, and integrating equation \eqref{Eqn_InfDomain_Heat} in time we see that $w_{\infty}(\bx) = \int_{0}^{\infty} w(t,\bx) \;dt$. However, the problem \eqref{Eqn_InfDomain_Heat} is posed over an infinite domain, and a truncation of the infinite domain $\mathbb{R}^d$ is inevitable for computational reasons. Therefore, we consider the problem 
\begin{align} \label{Eqn_FiniteDomain_Heat}
\partial_t v_{R}(t,\bx)-\Delta v_{R}(t,\bx)&= 0, \quad \text{in} \quad  (0,T) \times K_R:=[-R,R]^d, \\
v_{R}(0,\bx)  &= F(\bx). \nonumber \\
v_{R}(t,\bx)  &= 0 ´\quad \text{on} \quad (0,T) \times \partial K_R,  \nonumber
\end{align}
and define $v_{T,R}(\bx) := \int_{0}^{T} v_{R}(t,\bx) \;dt$. Then integrating equation \eqref{Eqn_FiniteDomain_Heat} in time, and using the fact that $v_{R}(T,\cdot) = e^{T  \Delta } F(\cdot)$, we obtain \eqref{Eqn_RegularisedElliptic}, which contains a regularization term in the form of exponential power of Laplace operator. In fact, this regularization term makes the boundary error (due to the artificial homogeneous Dirichlet boundary conditions in \eqref{Eqn_RegularisedElliptic}) decay at a Gaussian rate with increasing size, $R$, of the computational geometry. 
\begin{Remark}
    One immediately observes that upon choosing $T = \infty$ in \eqref{Eqn_RegularisedElliptic}, the exponential regularisation term disappears and one recovers a truncated domain version of the original PDE \eqref{Eqn_InfDomain_Elliptic}. This already indicates the fact that having a finite value for the parameter $T$ may be necessary for getting a rapidly decreasing boundary error.
\end{Remark}

\section{Analysis}
\label{Sec_Analysis}
In this section, we aim at analysing the following problems:
\begin{itemize}
    \item[a.] What is the error estimate when the infinite domain problem \eqref{Eqn_InfDomain_Elliptic} is replaced by a finite domain problem equipped with homogeneous Dirichlet boundary conditions ? 
    \item[b.] Error analysis of the proposed approach in the case of an infinitely large periodic magnet.  
    \item[c.] Error analysis of the proposed approach in the case of a right hand side $F$ with vanishing Fourier transform in the origin $\bome = {\bf 0}$.   
\end{itemize}

In particular, the analysis above will show that artificial homogeneous boundary conditions posed for the case $a$ will result in a large error in the region of interest, where the magnetic body exists. While the latter two cases will demonstrate the exponential decay of the artificial boundary conditions in the interior region under two different theoretical settings. Note that the last theoretical setting is very much of practical interest and is seen in a wide variety of applications within high-frequency magnetism, \cite{Jiale_etall_2016,XU2022169815}. 

\subsection{Effect of truncating the infinite domain in \eqref{Eqn_InfDomain_Elliptic}}
Let us consider the problem \eqref{Eqn_InfDomain_Elliptic} but posed over a finite domain, say $K_R:=[-R,R]^d$, $R\gg 1$ equipped with homogeneous Dirichlet boundary conditions.
\begin{align} \label{Eqn_FiniteDomain_Elliptic}
-\Delta w_{R}(\bx) &= F(\bx), \quad \bx \in K_R, \\
w_{R}(\bx) & = 0, \quad \bx \in \partial K_R  \nonumber.
\end{align}
We want to understand the decay of the error $w_{\infty}(\bx) - w_{R}(\bx)$ when such a simple truncation of $\mathbb{R}^d$ is assumed. Note that we are interested in the error only inside the interior region, say $K= \Omega =[-1/2,1/2]^d$, filled with the magnetic material. The following Theorem shows the poor decay ($\approx 1/R$) of the error in $L^2(K)$ (seen in the limit as $\alpha \to 1$).

\begin{Theorem} \label{Thm_Analysis_1}
    Let $F \in L^{\infty}(K)$ and assume that $w_{\infty}$ and $w_{R}$ solve the equations \eqref{Eqn_InfDomain_Elliptic} and \eqref{Eqn_FiniteDomain_Elliptic} respectively. Moreover, let $\alpha$ be any constant in the interval \ $0<\alpha <1$, and assume that $d \geq 3$. Then we have the following asymptotic error estimate (as $R \to \infty$) 
    \begin{align*}
            \| w_{\infty}- w_{R}(\bx) \|_{L^2(K)} \leq C R^{-d +3 - \alpha}  \| F \|_{L^{\infty}(K)},
    \end{align*}
    where $C$ is a constant which depends only on the dimension. 
\end{Theorem}
\begin{Proof} [Proof of Theorem \ref{Thm_Analysis_1}]
    We define a cutt-off function $\rho \in C_{c}^{\infty}(K_R)$ such that
    \begin{align*}
        \rho(\bx) = \begin{cases} 1 & \text{ if } \bx \in K_{R - \frac{1}{4} R^{\alpha}} \\
        0 & \text{ if }  \bx \text{ on } \partial K_{R}, 
        \end{cases}
    \end{align*}  
and $|\nabla  \rho | < \frac{C}{R^{\alpha}}$ on $E := K_{R} \backslash K_{R-\frac{1}{4}R^{\alpha}}$. Let us define $\tilde{e}_{R}(\bx) := \rho w_{\infty}(\bx) - w_{R}(\bx)$. Then clearly $\tilde{e}_{R}(\bx) = w_{\infty}(\bx) - w_{R}(\bx)$ in $K$ since $R\gg 1$. Moreover, 
\begin{align} \label{Eqn_Difference_Elliptic}
-\Delta \tilde{e}_{R}(\bx) &= - w_{\infty} \Delta \rho - 2 \nabla \rho \cdot \nabla w_{\infty} + (1-\rho) F, \quad \bx \in K_R, \\
\tilde{e}_{R}(\bx) & = 0, \quad \bx \in \partial K_R  \nonumber.
\end{align}
Now, let $G_{R}(\bx,\by)$ be the Green's function corresponding to \eqref{Eqn_Difference_Elliptic}, so that we write:
\begin{align*}
    \tilde{e}_{R}(\bx) &= \int_{K_R} G_{R}(\bx,\by) \left( - w_{\infty} \Delta \rho - 2 \nabla \rho \cdot \nabla w_{\infty} + (1-\rho) F \right)(\by) \; d\by \\
    &= \int_{K_R} w_{\infty}(\by) \nabla G_{R}(\bx,\by) \cdot \nabla \rho(\by)  - G_{R}(\bx,\by) \nabla \rho(\by) \cdot \nabla w_{\infty}(\by) \; d\by \\
    &= \int_{E} w_{\infty}(\by) \nabla G_{R}(\bx,\by) \cdot \nabla \rho(\by)  - G_{R}(\bx,\by) \nabla \rho(\by) \cdot \nabla w_{\infty}(\by) \; d\by.
\end{align*}
Note that in the second equality, we used the fact that $(1-\rho) F$ is zero on $E$ since $F$ is nonzero only on $K$. For all $\bx \in K$, we have 
\begin{align} \label{Estimate_L2Norm_StandardMethod}
    |\tilde{e}_{R}(\bx)| \leq  \dfrac{C}{R^{\alpha}}  \left( \|\nabla G_{R}(\bx,\cdot) \|_{L^{2}(E)} \| w_{\infty}(\cdot) \|_{L^{2}(E)} + |E| \|G_{R}(\bx,\cdot) \|_{L^{\infty}(E)} \| \nabla w_{\infty}(\cdot) \|_{L^{\infty}(E)}   \right)
\end{align}
Using a maximum principle of the form $0< G_{R}(\bx,\by) \leq G_{\infty}(\bx,\by) = C \frac{1}{|\bx - \by|^{d-2}}$, where $G_{\infty}$ is the Green's function corresponding to the infinite domain problem \eqref{Eqn_InfDomain_Elliptic}, we can bound
\begin{align*}
   \sup_{\bx \in K} \|G_{R}(\bx,\cdot) \|_{L^{\infty}(E)}  \leq \sup_{\bx \in K} \|G_{\infty}(\bx,\cdot) \|_{L^{\infty}(E)}  \leq C \dfrac{1}{|R-\frac{1}{4}R^{\alpha}-1|^{d-2}}
\end{align*}
Moreover, 
\begin{align} \label{Estimates_w_inf}
    \| w_{\infty}(\cdot) \|_{L^{\infty}(E)}  = \sup_{\bz \in E} \left| \int_{K} G_{\infty}(\bz,\by) F(\by) \; d\by \right| \leq C |K| \dfrac{1}{|R-\frac{1}{4}R^{\alpha}-1|^{d-2}} \| F \|_{L^{\infty}(K)}, 
\end{align}
\begin{align} \label{Estimates_w_inf_grad}
    \| \nabla w_{\infty}(\cdot) \|_{L^{\infty}(E)}  = \sup_{\bz \in E} \left| \int_{K} \nabla_{\bz} G_{\infty}(\bz,\by) F(\by) \; d\by \right| \leq C |K| \dfrac{1}{|R-\frac{1}{4}R^{\alpha}-1|^{d-1}} \| F \|_{L^{\infty}(K)}, 
\end{align}
It remains to bound the term $\|\nabla G_{R}(\bx,\cdot) \|_{L^{2}(E)}$. For this, we use a Cacioppoli type inequality. In particular, we prove that 
\begin{align} \label{Ineq_Cacciopoli_Elliptic}
    \int_{E} |\nabla_{\by} G_{R}(\bx,\by)|^2 \; d\by \leq \frac{C}{R^{2\alpha}} \int_{K_{R} \backslash K_{R-\frac{3R^{\alpha}}{4}}} | G_{R}(\bx,\by)|^2 \; d\by.
\end{align}
To do this, we define $\eta(\bx) \geq 0$ for all $x \in K_{R}$ as
\begin{align*}
    \eta(\bx) = \begin{cases} 1 & \quad \text{ in } \quad K_R \backslash K_{R - \frac{1}{4}R^{\alpha}} \\
    0 & \quad \text{ in } \quad K_{R-\frac{3R^{\alpha}}{4}}, 
    \end{cases}
\end{align*}
such that $|\nabla \eta| \leq \frac{C}{R^{\alpha}}$.
Next, we multiply the equation $-\Delta G_R(\bx,\by) = \delta(\bx-\by)$, equipped homogeneous boundary conditions on $\partial K_R$,  with the term $\eta^2 G_R$ and integrate  over $K_R$ to obtain 
\begin{align*}
    \int_{K_{R}} \nabla G_{R}(\bx,\by) \cdot \nabla (\eta^2(\by) G_{R}(\bx,\by)) \; d\by = 0.
\end{align*}
The right hand side is zero since $\bx \in K$ and $R>3$. We observe that 
\begin{align*}
   0= \int_{K_{R}} \nabla G_{R}(\bx,\by) \cdot \nabla (\eta^2(\by) G_{R}(\bx,\by)) \; d\by &= \int_{K_{R}} \nabla (\eta G_{R}) \cdot \nabla (\eta G_{R}) \; d\by \\ &- \int_{K_{R}} G_R^2 \nabla \eta  \cdot \nabla \eta \; d\by.
\end{align*}
Finally, using the fact that $|\nabla \eta| \leq \frac{C}{R^{\alpha}}$, it follows for all $x\in K$ that 
\begin{align*}
    \int_{E} |\nabla_{\by} G_{R}(\bx,\by)|^2 \; d\by & \leq \int_{K_{R}} \nabla (\eta G_{R}) \cdot \nabla (\eta G_{R}) \; d\by \leq \frac{C}{R^{2\alpha}} \int_{K_{R} \backslash K_{R-\frac{3}{4}R^{\alpha}}} | G_{R}(\bx,\by)|^2 \; d\by \\
    &\leq \dfrac{C}{R^{2\alpha}} R^{d-1+\alpha} \| G_{R}(\bx,\cdot) \|_{L^{\infty}(K_{R} \backslash K_{R-\frac{3}{4}R^{\alpha}})}^2 \\
    &\leq R^{d-1-\alpha} \dfrac{1}{|R-\frac{3}{4}R^{\alpha}-1|^{2d-4}}
\end{align*}
Using the inequalities \eqref{Estimates_w_inf}, \eqref{Estimates_w_inf_grad}, and \eqref{Ineq_Cacciopoli_Elliptic}, in \eqref{Estimate_L2Norm_StandardMethod}, and using the fact that $\| f \|_{L^2(E)} \leq |E|^{1/2} \| f \|_{L^{\infty}(E)}$, and $|E| = O(R^{d-1/2})$, together with the maximum principle $0< G_{R}(\bx,\by) \leq G_{\infty}(\bx,\by) = C \frac{1}{|\bx - \by|^{d-2}}$, we conclude that  
\begin{align*} 
    |e_{R}(\bx)| &\leq  \dfrac{C}{R^{\alpha}}  \left( R^{(d-1-\alpha)/2} \dfrac{1}{|R-\frac{3}{4}R^{\alpha}-1|^{d-2}} R^{(d-1-\alpha)/2} \dfrac{1}{|R-\frac{1}{4}R^{\alpha}-1|^{d-2}} \| F \|_{L^{\infty}(K)} \right)  \\
    &+ \dfrac{C}{R^{\alpha}}  \left( R^{d-1-\alpha} \dfrac{1}{|R-\frac{1}{4}R^{\alpha}-1|^{d-2}}  \dfrac{1}{|R-\frac{1}{4}R^{\alpha}-1|^{d-1}} \| F \|_{L^{\infty}(K)}   \right).
\end{align*}
Since $0<\alpha<1 $, and $R \gg 1$, we get the following overall asymptotic rate:
\begin{align*} 
    |e_{R}(\bx)| &\leq \dfrac{C}{R^{\alpha}}  R^{d-1} \dfrac{1}{R^{2d-4}} \| F \|_{L^{\infty}(K)} + \dfrac{C}{R^{d-2}}  \| F \|_{L^{\infty}(K)} \\
    &\leq C R^{-d +3 - \alpha}  \| F \|_{L^{\infty}(K)}.
\end{align*}

$\square$

\end{Proof}

\subsection{Analysis of an infinitely large periodic magnet}
Here we  assume that $\bM \in W^{1}_{per}(K)$ is a periodic function, say with period $1$ in all directions, so that $F = - \nabla \cdot \bM(\bx)$ for all $\bx \in \mathbb{R}^d$. Therefore, the solution $u(\bx)$ in \eqref{Eqn_Hybrid_Representation} is given only by $v_{T,R}$ and the boundary integral term is no longer needed. Here, we intend to give an error estimate for the difference between the solution, $v_{T,R}$, to the equation \eqref{Eqn_RegularisedElliptic}, and the solution $w_{\infty}$ to equation \eqref{Eqn_InfDomain_Elliptic} for a periodic right hand side $F$. Due to the periodicity of the problem, we assume that the region of interest is the unit cube $K:= [-1/2,1/2]^d$ in $\mathbb{R}^d$, and we have to show that the difference $w_{\infty}(\bx) - v_{T,R}(\bx)$ is small for $\bx \in K$. To be able to do this, we split the error $w_{\infty} - v_{T,R}$ as follows: 

\begin{align}\label{Eqn_Error_Splitting}
    w_{\infty}(\bx) - v_{T,R}(\bx)  = \underbrace{\left( w_{\infty}(\bx) - w_T(\bx) \right)}_{\text{Modelling error}}  + \underbrace{\left(    w_T(\bx) - v_{T,R}(\bx) \right)}_{\text{Boundary error}},
\end{align}
where $w_T(\bx) = \int_{0}^{T} w(t,\bx) \; dt$, and $w$ solves \eqref{Eqn_InfDomain_Heat}. The first term is the modelling error due to the integration of $w$ over a finite time $T$. As explained earlier, the term $e^{T\Delta} F$ in equation \eqref{Eqn_RegularisedElliptic} is due to this finite time integration, which act as an extra regularisation term in comparison to the original infinite domain model problem \eqref{Eqn_InfDomain_Elliptic}; hence the name modelling error. The second term, on the other hand, is the boundary error due to a truncation of the infinite domain $\mathbb{R}^d$ to $K_{R}$. 

The main result of this section is the following Theorem, which shows the exponential decay of the overall error as the domain size $|K_{R}|$ increases. 

\begin{Theorem} \label{Thm_Analysis_b}
    Assume that $F \in L^2(K)$ and that $T \leq \frac{(R-1)^2}{2}$. Then we have 
    \begin{align*}
        \| w_{\infty} - v_{T,R} \|_{L^2(K)} \leq C_1 \left( \dfrac{1}{\lambda_1} e^{-\lambda_1 T} +  R^{d-1} \left( 1 + C_2 \dfrac{T}{|R-1|} \right) \dfrac{1}{T^{d/2}} e^{-\frac{(R-1)^2}{4T}}  \right) \|  F  \|_{L^2(K)},
    \end{align*}
    where $C_{1,2}$ do not depend on $R$ and $T$ but may depend on the dimension. Moreover, upon choosing $T=\frac{R-1}{2\sqrt{\lambda_1}}$, we have
 \begin{align*}
        \| w_{\infty} - v_{T,R} \|_{L^2(K)} \leq C   R^{d-1} \dfrac{1}{T^{d/2}} e^{-\frac{\sqrt{\lambda_1} (R-1)}{2}} \|  F  \|_{L^2(K)},
    \end{align*}
    where $C$ is a constant independent of $R$ and $T$.
\end{Theorem}

\begin{Proof} [Proof of Theorem \ref{Thm_Analysis_b}]The proof of the above theorem follows immediately from  Lemmas \ref{Lemma_Per_ModError} and \ref{Lemma_Per_BoundError}, which are given to bound the modelling and the boundary errors separately. 

$\square$
\end{Proof}
\begin{Lemma} (Modelling error) \label{Lemma_Per_ModError}
    Assume that $F \in L^2(K)$. Then we have 
    \begin{align*}
        \| w_{\infty} - w_T \|_{L^2(K)} \leq \dfrac{1}{\lambda_1} e^{-\lambda_1 T} \| F \|_{L^2(K)},
    \end{align*}
    where $\lambda_1 >0$ is the smallest positive eigenvalue of the operator $-\Delta$ equipped with periodic boundary conditions on $K$. 
\end{Lemma}
\begin{Proof} [Proof of Lemma \ref{Lemma_Per_ModError}] We write 
\begin{align*}
    \| w_{\infty} - w_T \|_{L^2(K)}  &= \| \int_{0}^{\infty} w(t,\cdot) \; dt -  \int_{0}^{T} w(t,\cdot) \; dt\|_{L^2(\Omega)} \\
    &=\| \int_{T}^{\infty} w(t,\cdot) \; dt \|_{L^2(K)} \leq \int_{T}^{\infty} \|  w(t,\cdot) \|_{L^2(K)} \; dt.
\end{align*}

 Now we bound $\| w(t,\cdot) \|_{L^2(K)}$ as follows. Let $w(t,x) = \sum_{j=0}^{\infty} w_j(t) \varphi_{j}(\bx)$, and $F(\bx) = \sum_{j=0}^{\infty} F_j \varphi_j(\bx)$. Plugging the eigenfunction expansion of $w$ into the PDE, and using the orthogonality of eigenfunctions, we immediately see that 
    \begin{align*}
        w_j^{\prime}(t) + \lambda_j w_j(t) &= 0, \quad j=0,1,\ldots\\
        w_j(0) &= F_j, 
    \end{align*}
    where $w_0(0) = F_0 = 0$ due to $\langle F,1 \rangle_K = 0$, and 
    $w_j(t)  = e^{-\lambda_j t} F_j, \quad j\geq 1$. Therefore, 
     \begin{align*}
       \| w(t,\cdot) \|_{L^2(K)}^2  &= \sum_{j=1}^{\infty} w_j^2(t)  = \sum_{j=1}^{\infty} e^{-2\lambda_j t} |F_j|^2 \\
       &\leq e^{-2\lambda_1 t} \sum_{j=1}^{\infty} |F_j|^2   =  e^{-2\lambda_1 t} \|  F \|_{L^2(K)}^2, 
    \end{align*}
    where we used the Parseval's equality in the first and last steps. We obtain $\| w(t,\cdot) \|_{L^2(K)} \leq e^{-\lambda_1 t} \|  F \|_{L^2(K)}$ by taking the square root of both sides. Therefore,
\begin{align*}
    \| w_{\infty} - w_T \|_{L^2(K)}  &\leq \int_{T}^{\infty} \|  w(t,\cdot) \|_{L^2(K)} \; dt 
    \leq \int_{T}^{\infty} e^{-\lambda_1 t} \; dt  \| F \|_{L^2(K)}  \leq \dfrac{1}{\lambda_1} e^{-\lambda_1 T} \| F \|_{L^2(K)}. 
\end{align*}
    $\square$
\end{Proof}
\begin{Lemma} (Boundary error) \label{Lemma_Per_BoundError}
    Assume that $F \in L^2(K)$ and $T \leq \frac{(R-1)^2}{2}$. Then we have 
    \begin{align*}
        \| w_{T} - v_{T,R} \|_{L^2(K)} \leq C_1 \| F \|_{L^{2}(K)} R^{d-1} \left( 1 + C_2 \dfrac{T}{|R-1|} \right) \dfrac{1}{T^{d/2}} e^{-\frac{(R-1)^2}{4T}},
    \end{align*}
    where $C_1$ and $C_2$ are constants independent of $R$ and $T$ but may depend on $d$.
\end{Lemma}
\begin{Proof}[Proof of Lemma \ref{Lemma_Per_BoundError}]
    We define $e_{T,R}(\bx)  = w_T(\bx)  - v_{T,R}(\bx) = \int_{0}^{T} w(t,\bx) - v_{R}(t,\bx)\; dt$, where $v_R$ solves \eqref{Eqn_FiniteDomain_Heat}, and observe that $e_{R}(\bx) = \int_{0}^{T} e_{R}(t,\bx) \; dt$, where 
    \begin{align*}
        \partial_t e_{R}(t,\bx) - \Delta e_{R}(t,\bx) &= 0, \quad \text{in} \quad  (0,T) \times K_R\\
        e_R(0,\bx) &= 0, \quad \text{in} \quad  K_R\\
        e_R(t,\bx) &= w(t,\bx) \quad \text{on} \quad (0,T) \times \partial K_R.
    \end{align*}
    Let $\te_{R}(t,\bx) = \rho w(t,\bx) - v_{R}(t,\bx)$, where $\rho \in C_{c}^{\infty}(K_R)$ is a cutt-off function which satisfies 
    \begin{align*}
        \rho(\bx) = \begin{cases} 1 & \text{ if } \bx \in K_{\tilde{R}} \\
        0 & \text{ if }  \bx \text{ on } \partial K_{R}, 
        \end{cases}
    \end{align*}    
    where $R-1<\tilde{R} < R-1/2$. Then
    \begin{align*}
        \| e_{T,R} \|_{L^2(K)}  &=  \| \int_{0}^{T} w(t,\cdot) - v_{R}(t,\cdot)  \; dt \|_{L^2(K)} 
        \leq \int_{0}^{T} \| w(t,\cdot) - v_{R}(t,\cdot) \|_{L^2(K)} \; dt \\ &=  \int_{0}^{T} \| \rho w(t,\cdot) - v_{R}(t,\cdot) \|_{L^2(K)} \; dt  = \int_{0}^{T} \| \te_{R}(t,\cdot) \|_{L^2(K)} \; dt \\
        &\leq |K| \int_{0}^{T} \| \te_{R}(t,\cdot) \|_{L^{\infty}(K)} \; dt 
    \end{align*}
    Moreover, we use the following estimate from \cite{Abdulle_Arjmand_Paganoni_18b} (Lemma 4.5 in \cite{Abdulle_Arjmand_Paganoni_18b}):
     \begin{align*}
        \| \te_{R}(t,\cdot) \|_{L^{\infty}(K)} \leq C_1 \| F \|_{L^{2}(K)} \left( 1 + C_2 \dfrac{t}{|R-1|} \right) \dfrac{R^{d-1}}{t^{d/2}} e^{-\frac{(R-1)^2}{4t}}, \quad \text{ whenever } t < (R-1)^2
    \end{align*}
    to conclude that 
    \begin{align*}
        \| e_{T,R} \|_{L^2(K)}  &\leq |K| \int_{0}^{T} \| \te_{R}(t,\cdot) \|_{L^{\infty}(K)} \; dt  \\
        &\leq C_1 \| F \|_{L^{2}(K)} R^{d-1} \left( 1 + C_2 \dfrac{T}{|R-1|} \right) \int_{0}^{T} \dfrac{1}{t^{d/2}} e^{-\frac{(R-1)^2}{4t}} \; dt \\
        &\leq C_1 \| F \|_{L^{2}(K)} R^{d-1} \left( 1 + C_2 \dfrac{T}{|R-1|} \right) \dfrac{1}{T^{d/2}} e^{-\frac{(R-1)^2}{4T}}.
    \end{align*}
    Note that in the last estimate, we used the fact that $\dfrac{1}{t^{d/2}} e^{-\frac{(R-1)^2}{4t}}$ is an increasing function of $t$ whenever $t < \frac{|R-1|^2}{2}$. 
$\square$
\end{Proof}

\subsection{Analysis for magnetization with vanishing frequency components in the origin}

Here, we relax the somewhat theoretical periodicity assumption in the previous section and we instead assume that $F \in \mathcal{H}_{\omega_0}$, and we provide error estimates both for the modelling and the boundary errors. Such an assumption on $F$ is much more realistic due to the oscillatory nature of the magnetization vector in the micromagnetics regime. Note that the analysis in this section relies on the splitting \eqref{Eqn_Error_Splitting}, and the terms $w_{\infty},w_T,v_{T,R}$ solve the same equations as in the previous section with the exception that now $F = -\nabla \cdot \bM$ is non-zero only over $K=[-1/2,1/2]^d$. Therefore, similar to the periodic analysis we are interested in analyzing the error over $K$.

The main result of this section is the following Theorem, which again shows an exponentially decaying errors in $R$ upon optimally choosing the parameter $T$.

\begin{Theorem} \label{Thm_Main_Analysis_c}
    Assume that $F \in \mathcal{H}_{\omega_0}$ and $\nabla F \in L^{2}(\mathbb{R}^d)$. Moreover, let $T < |R-1|^2$. Then we have 
    \begin{align*}
        \| w_{\infty} - v_{T,R} \|_{L^2(K)} \leq \dfrac{1}{\omega_0^2} e^{-\omega_0^2 T} \| F \|_{L^2(\mathbb{R}^{d})} + C_{R,T,\omega_0}  e^{-c_d \frac{|R-1|^2}{T} } \left( \| F \|_{L^2(\mathbb{R}^d)}  + \| \nabla F \|_{L^2(\mathbb{R}^d)}\right),
    \end{align*} 
    where $c_d$ may depend only on the dimension $d$, and $C_{R,T,\omega_0}  = C \dfrac{R^{(d-1)/2}}{T^{d/2}} \left( T + T^2 \right) \left( \frac{1}{|R-1|\omega_0^2} + 1 \right)$, for a constant $C$ independent of $R,T,\omega_0$ and $F$. Moreover, the optimal value of $T = \frac{|R-1|\sqrt{c_d}}{\omega_0}$, for which the estimate simplifies to
      \begin{align*}
        \| w_{\infty} - v_{T,R} \|_{L^2(K)} \leq \left( \dfrac{1}{\omega_0^2} + C_{R,T,\omega_0} \right) e^{-\omega_0 \sqrt{c_d} |R-1|} \left( \| F \|_{L^2(\mathbb{R}^d)}  + \| \nabla F \|_{L^2(\mathbb{R}^d)}\right).
    \end{align*} 
\end{Theorem}
\begin{Proof}[Proof of Theorem \ref{Thm_Main_Analysis_c}] The proof of this Theorem is a direct consequence of the error splitting \eqref{Eqn_Error_Splitting} and Lemmas \ref{Lemma_Modelling_Error_Analysis_c} and \ref{Lemma_Boundary_Error_Analysis_c} below, where separate estimates for the  modelling and the boundary errors are given. 
\end{Proof}

\begin{Lemma} \label{Lemma_Estimate_w_Freq}
    Suppose that $w$ solves the following heat equation 
\begin{align*} 
\partial_t w(t,\bx)-\Delta w(t,\bx)&= 0, \quad \bx \in \mathbb{R}^d, \\
w(0,\bx)  &= F(\bx), \nonumber
\end{align*}
where $F \in \mathcal{H}_{\omega_0}$. Then the following estimate holds:
\begin{align*}
    \|  w(t,\cdot) \|_{L^2(\mathbb{R}^d)} \leq e^{-\omega_0^2 t} \|  F  \|_{L^2(\mathbb{R}^d)}.
\end{align*}
Moreover, if in addition $\nabla F \in L^2(\mathbb{R}^d)$, then 
\begin{align*}
    \|  \nabla w(t,\cdot) \|_{L^2(\mathbb{R}^d)} \leq e^{-\omega_0^2 t} \|  \nabla F  \|_{L^2(\mathbb{R}^d)}.
\end{align*}
\end{Lemma}
\begin{Proof}[Proof of Lemma \ref{Lemma_Estimate_w_Freq}]
We start by expressing the solution $w$ in terms of the Green's function
\begin{align*}
        w(t,\bx) = \left( G(t,\bx,\cdot) \star F(\cdot) \right)(\bx).
\end{align*}
Therefore, 
\begin{align*}
        \| w(t,\cdot) \|_{L^2(\mathbb{R}^d)} &= \| \hat{G} \hat{F} \|_{L^2(\mathbb{R}^d)} = \sqrt{\int_{\mathbb{R}^d} e^{-2|\bome|^2 t} |\hat{F}(\bome)|^2 \; d\bome } \\
        &= \sqrt{\int_{|\bome| > \omega_0 } e^{-2|\bome|^2 t} |\hat{F}(\bome)|^2 \; d\bome } \\
        &\leq e^{-|\omega_0|^2 t} \sqrt{\int_{|\bome| > \omega_0 } |\hat{F}(\bome)|^2 \; d\bome } \leq e^{-\omega_0^2 t} \| F \|_{L^2(\mathbb{R}^d)}.
\end{align*}
To prove the second statement we proceed similarly 
\begin{align*}
        \| \nabla w(t,\cdot) \|_{L^2(\mathbb{R}^d)} &= \| \widehat{\nabla G} \hat{F} \|_{L^2(\mathbb{R}^d)} = \sqrt{\int_{\mathbb{R}^d} |\bome|^2 e^{-2|\bome|^2 t} |\hat{F}(\bome)|^2 \; d\bome } \\
        &= \sqrt{\int_{|\bome| > \omega_0 } e^{-2|\bome|^2 t} |\bome \hat{F}(\bome)|^2 \; d\bome } \\
        &\leq e^{-|\omega_0|^2 t} \sqrt{\int_{|\bome| > \omega_0 } |\bome \hat{F}(\bome)|^2 \; d\bome } \leq e^{-\omega_0^2 t} \| \nabla F \|_{L^2(\mathbb{R}^d)}.
\end{align*}
\end{Proof}
\begin{Lemma} \label{Lem_Intermediate}
    Suppose that $G_p$ is the Green's function of the parabolic PDE \eqref{Eqn_FiniteDomain_Heat}, then 

\begin{align}
    \| G_p(t,x;\cdot)\|_{L^2(E)} \leq C \dfrac{|E|^{1/2}}{t^{d/2}} e^{-c_d \frac{|R-1|^2}{t}}
\end{align}
  and 
  \begin{align}
    \| G_p(\cdot,x;\cdot)\|_{L^2((0,t)\times E)} \leq C \dfrac{|E|^{1/2}}{|R-1|t^{d/2-1}} e^{-c_d \frac{|R-1|^2}{t}},
\end{align}
where $c_d$ and $C$ are constants depending only on the dimension $d$. 
\end{Lemma}
\begin{Proof}[Proof of Lemma \ref{Lem_Intermediate}] This Lemma is a natural consequence of the Nash-Aronson estimate, \cite{Aronson_67}, and a complete proof can be found in \cite{Abdulle_Arjmand_Paganoni_18b}. 
\end{Proof}

\begin{Lemma} \label{Lemma_Modelling_Error_Analysis_c} (Modelling error)
    Assume that $F \in \mathcal{H}_{\bome_0}$. Then we have 
    \begin{align*}
        \| w_{\infty} - w_T \|_{L^2(K)} \leq \dfrac{1}{\omega_0^2} e^{-\omega_0^2 T} \| F \|_{L^2(\mathbb{R}^{d})}.
    \end{align*} 
\end{Lemma}
\begin{Proof}[Proof of Lemma \ref{Lemma_Modelling_Error_Analysis_c}] We write 
\begin{align*} 
    \| w_{\infty} - w_T \|_{L^2(K)}  &= \| \int_{0}^{\infty} w(t,\cdot) \; dt -  \int_{0}^{T} w(t,\cdot) \; dt\|_{L^2(K)} \\
    &=\| \int_{T}^{\infty} w(t,\cdot) \; dt \|_{L^2(K)} \leq \int_{T}^{\infty} \|  w(t,\cdot) \|_{L^2(\mathbb{R}^{d})} \; dt \\
    &\leq \int_{T}^{\infty} e^{-\omega_0^2 t} \; dt  \| F \|_{L^2(\mathbb{R}^{d})}  \leq \dfrac{1}{\omega_0^2} e^{-\omega_0^2 T} \| F \|_{L^2(\mathbb{R}^{d})}. 
\end{align*}
Note that we used Lemma \ref{Lemma_Estimate_w_Freq} to bound $w$ in $L^2$.
    $\square$
\end{Proof}
\begin{Lemma} (Boundary error) \label{Lemma_Boundary_Error_Analysis_c}
    Assume that $F \in \mathcal{H}_{\omega_0}$ and $\nabla F \in L^{2}(\mathbb{R}^d)$. Moreover, let $t < |R-1|^2$. Then we have 
    \begin{align*}
        \| w_{T} - v_{T,R} \|_{L^2(K)} \leq C_{R,T,\omega_0}  e^{-c_d \frac{|R-1|^2}{T} } \left( \| F \|_{L^2(\mathbb{R}^d)}  + \| \nabla F \|_{L^2(\mathbb{R}^d)}\right).
    \end{align*}
    where $C_{R,T,\omega_0}  = C \dfrac{R^{(d-1)/2}}{T^{d/2}} \left( T + T^2 \right) \left( \frac{1}{|R-1|\omega_0^2} + 1 \right)$, for a constant $C$ independent of $R,T,\omega_0$ and $F$. 
\end{Lemma}
\begin{Proof}
    The boundary error is given by the difference $e_{T,R}(\bx)  = w_T(\bx)  - v_{T,R}(\bx) = \int_{0}^{T} w(t,\bx) - v_{R}(t,\bx)\; dt$. Note that once again $e_{T,R}(\bx) = \int_{0}^{T} e_{R}(t,\bx) \; dt$, where 
    \begin{align*}
        \partial_t e_{R}(t,\bx) - \Delta e_{R}(t,\bx) &= 0, \quad \text{in} \quad  (0,T) \times K_R\\
        e_R(0,\bx) &= 0, \quad \text{in} \quad  K_R\\
        e_R(t,\bx) &= w(t,\bx) \quad \text{on} \quad (0,T) \times \partial K_R.
    \end{align*}
    Similar to the analysis in the previous section, let $\te_{R}(t,\bx) = \rho w(t,\bx) - v_{R}(t,\bx)$, where $\rho \in C_{c}^{\infty}(K_R)$ is a cutt-off function which satisfies 
    \begin{align*}
        \rho(\bx) = \begin{cases} 1 & \text{ if } \bx \in K_{\tilde{R}} \\
        0 & \text{ if }  \bx \in \partial K_{R}, 
        \end{cases}
    \end{align*}    
    and $|\partial^2_{x_j}   \rho | < C_{\rho}$ in $E:=K_{R}\backslash K_{\tilde{R}}$, for all $j=1,\ldots,d$. Then
     \begin{align*}
        \partial_t \te_{R}(t,\bx) - \Delta \te_{R}(t,\bx) &= - 2 \nabla \rho \cdot \nabla w - w \Delta \rho, \quad \text{in} \quad  (0,T) \times K_R\\
        \te_R(0,\bx) &= (\rho - 1)F, \quad \text{in} \quad  K_R\\
        \te_R(t,\bx) &= 0 \quad \text{on} \quad (0,T) \times \partial K_R.
    \end{align*}
    Therefore, with $f(s,\by) = - 2 \nabla \rho(\by) \cdot \nabla w(s,\by) - w(s,\by) \Delta \rho(\by)$, and for all $\bx \in K$ we have
      \begin{align*}
        \left| \te_{R}(t,\bx)  \right| &= \left| \int_{K_R} \int_{0}^{t} G_{p}(t,\bx;s,\by) f(s,\by) \; ds d\by + \int_{K_R} G_{p}(t,\bx;0,\by) (\rho - 1)F(\by) \; d\by \right| \\
        &= \left| \int_{E} \int_{0}^{t} G_{p}(t,\bx;s,\by) f(s,\by) \; ds d\by + \int_{E} G_{p}(t,\bx;0,\by) (\rho - 1)F(\by) \; d\by \right| \\
        &\leq C \| G_p(t,\bx;\cdot,\cdot) \|_{L^2((0,t)\times E)} \left(\| \nabla w \|_{L^2((0,t)\times E)}  + \| w \|_{{L^2((0,t)\times E)}}\right) + \| G_p(t,\bx;0,\cdot) \|_{L^2(E)} \| F\|_{L^2(E)} \\
        &\leq C \dfrac{|E|^{1/2}}{|R-1|t^{d/2-1}} e^{-c_d \frac{|R-1|^2}{t}} \left( \| \nabla w \|_{L^2((0,t)\times \mathbb{R}^d)} +  \| w \|_{L^2((0,t)\times \mathbb{R}^d)} \right) + C \dfrac{|E|^{1/2}}{ t^{d/2}} e^{-c_d \frac{|R-1|^2}{t}} \| F \|_{L^2(\mathbb{R}^d)} \\
        &\leq C \dfrac{|E|^{1/2} (1-e^{-\omega_0^2 t})}{|R-1| w_0^2 t^{d/2-1}} e^{-c_d \frac{|R-1|^2}{t}} \left( \| F\|_{L^2(\mathbb{R}^d)} +  \| \nabla F \|_{L^2(\mathbb{R}^d)} \right) \\ &+ C \dfrac{|E|^{1/2}}{ t^{d/2}} e^{-c_d \frac{|R-1|^2}{t}} \| F \|_{L^2(\mathbb{R}^d)} \\
        &\leq C \dfrac{|E|^{1/2}}{|R-1|w_0^2 t^{d/2-1}} e^{-c_d \frac{|R-1|^2}{t}} \left( \| F\|_{L^2(\mathbb{R}^d)} +  \| \nabla F \|_{L^2(\mathbb{R}^d)} \right) + C \dfrac{|E|^{1/2}}{ t^{d/2}} e^{-c_d \frac{|R-1|^2}{t}} \| F \|_{L^2(\mathbb{R}^d)}
    \end{align*}
    Therefore, it follows that 
    \begin{align*}
        \| e_{T,R} \|_{L^2(K)}  &=  \| \int_{0}^{T} w(t,\cdot) - v_{R}(t,\cdot)  \; dt \|_{L^2(K)} 
        \leq \int_{0}^{T} \| w(t,\cdot) - v_{R}(t,\cdot) \|_{L^2(K)} \; dt \\ &=  \int_{0}^{T} \| \rho w(t,\cdot) - v_{R}(t,\cdot) \|_{L^2(K)} \; dt  = \int_{0}^{T} \| \te_{R}(t,\cdot) \|_{L^2(K)} \; dt \\
        &\leq |K| \int_{0}^{T} \| \te_{R}(t,\cdot) \|_{L^{\infty}(K)} \; dt  \\
        &\leq C \dfrac{R^{(d-1)/2}}{T^{d/2}} \left( T + T^2 \right) \left( \frac{1}{|R-1|\omega_0^2} + 1 \right) e^{-c_d \frac{|R-1|^2}{T} } \left( \| F \|_{L^2(\mathbb{R}^d)}  + \| \nabla F \|_{L^2(\mathbb{R}^d)}\right).
    \end{align*}
$\square$
\end{Proof}

\section{Numerical Approximations}
\label{Sec_Num_Res}
Numerical approximation of the full algorithm consists of two parts: solving the regularised PDE \eqref{Eqn_RegularisedElliptic} and approximation of the PDE \eqref{Eqn_b} whose boundary values are given a boundary integral. Solving elliptic PDEs by far is standard and therefore we skip discussions in this direction. However, interesting is the numerical approximation of the exponential power of the Laplacian  in \eqref{Eqn_RegularisedElliptic} as well as that of the boundary integral \eqref{Eqn_b}. These issues are discussed in the upcoming subsections.
\subsection{Approximation of the regularised PDE \eqref{Eqn_RegularisedElliptic}}
We build our discussion based on the assumption that the PDE \eqref{Eqn_RegularisedElliptic} is discretized by a finite difference approximation. Such an approximation results in the following system to be solved: 
\begin{align*}
    -\Delta_{h} v_{h}  = F_h - e^{T \Delta_h} F_h.
\end{align*}
It is assumed that $\Delta_h \in \mathbb{R}^{N^d \times N^d}$. Central for efficiently solving this system is the approximation of the exponential term, which using a naive approximation can be extremely costly. To approximate this term, we use instead a low rank approximation based on the Arnoldi algorithm which reads as
\begin{equation*}
    e^{\mathcal{A}_h} F_h  \approx Q e^{\mathcal{H}} Q^{\star} F_h, \text{ where } \mathcal{A}_h  = -T \Delta_h,
\end{equation*}
where $Q\in \mathbb{R}^{N^d\times k}$, $\mathcal{H} \in \mathbb{R}^{k\times k}$ is an upper-Hessenberg matrix with $k \ll N^d$. For elliptic operators (such as $-\Delta_{h}$), the optimal value of $k$ is given by $k = O(N)$, while still retaining exponentially decaying bounds in terms of $R$, see \cite{Abdulle_Arjmand_Paganoni_2023} for more details. This implies that the value $k$, characterising the size of the low-rank matrix $\mathcal{H}$ is independent of the dimension and the exponential term can be accurately approximated at the cost of computing the exponential power of a matrix in $1$-dimension. In particular, the cost of the Arnoldi algorithm itself for this specific setting is $O(N^d k)$ implying an overall cost of $O(N^{d+1})$, which is optimal in three dimensions in the sense that the cost will be equal to the cost of approximating the boundary integral in \eqref{Eqn_b}, which scales as $O(N^{2d-2})$. 

\subsection{Numerical approximation of the boundary integral}
\begin{figure}[ht]
    \begin{center}
        \includegraphics[scale=0.4]{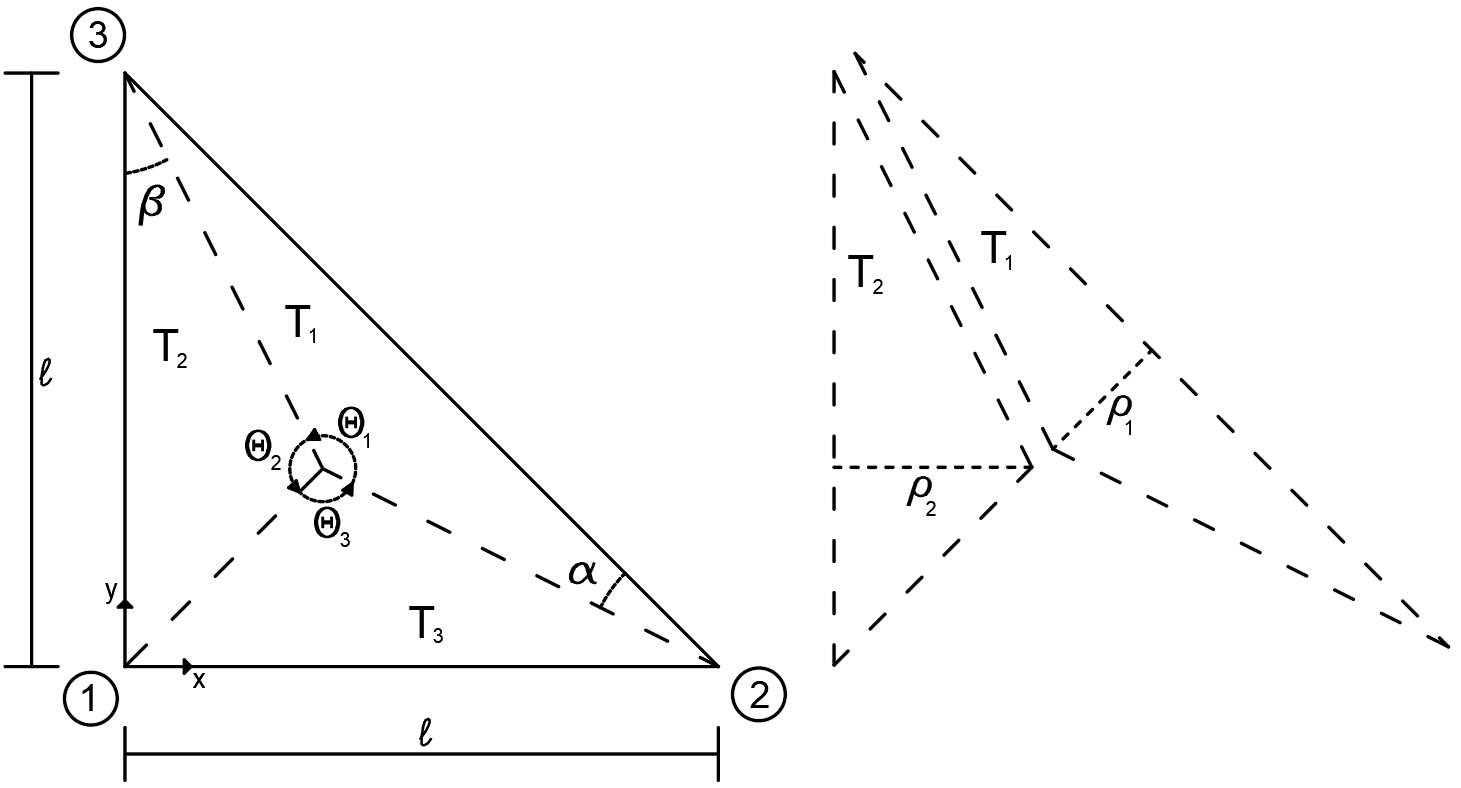}
        \caption{Transformation to polar coordinates in reference triangle.}
        \label{triangle}
    \end{center}
\end{figure}

Green functions for 2D and 3D Laplace's operator are
\begin{equation*}
    G(\bx,\bxp) = -\frac{1}{2\pi}\ln{|\bx-\bxp|} \qquad\text{and}\qquad G(\bx,\bxp) = \frac{1}{4\pi}\frac{1}{|\bx-\bxp|}
\end{equation*}
respectively, and they are used in \eqref{Eqn_Formal_Formulation}. The integral
\begin{equation}
    \int_{\partial\Omega} G(\bx,\bxp) \bM(\bxp) \cdot \bnp \; d\bxp
    \label{int_green}
\end{equation}
in three dimensions is calculated using triangles over the surface. The triangles used here have flat shapes defined by three vertex nodes. These triangles are in a parametric form and maps the triangle in the physical three-dimensional space to a right isosceles triangle in the $zw$ parametric plane. The mapping function from physical to parametric space is
\begin{equation*}
    \bx = \bx_1\gamma + \bx_2z + \bx_3w
\end{equation*}
where $\gamma=1-z-w$. Tangential vectors in the $zw$ plane, and their corresponding lengths are described, respectively, by
\begin{eqnarray*}
    \be_z & = & \bx_2-\bx_1 \qquad \be_w = \bx_3-\bx_1 \\
    h_z & = & |\bx_2-\bx_1| \qquad h_w = |\bx_3-\bx_1|
\end{eqnarray*}
The normal vector to $zw$ plane is
\begin{equation*}
    \bn = \frac{\be_z\times\be_w}{|\be_z\times\be_w|}    
\end{equation*}
The integral of the function $f(\bx,\bxp)=G(\bx,\bxp)\bM(\bxp)\cdot\bnp$ over the surface of the triangle in physical space is given by
\begin{equation*}
    \int_\triangle f(\bx,\bxp)d\bxp=|\be_z\times\be_w|\int_0^1\int_0^{1-z}f(\bx(z,w),\bxp(z,w))dwdz
\end{equation*}
When the surface function $f(\bx,\bxp)$ is non-singular over the surface of the triangle, the last integral could be computed accurately using a Gaussian quadrature rule for the surface of a triangle, for instance
\begin{equation*}
    \int_\triangle f(\bx,\bxp)d\bxp\approx\frac{|\be_z\times\be_w|}{2}\sum_{i=1}^{NGP}f(\bx(z_i,w_i),\bxp(z_i,w_i))\omega_i
\end{equation*}
where $NGP$ is the number of Gaussian quadrature points, the pair $(z_i,w_i)$ are the coordinates of the $i-th$ Gaussian point inside the triangle in the parametric space, and $\omega_i$ are the integration weights corresponding to the $i-th$ Gaussian point. \\
On the other hand, when surface function $f(\bx,\bxp)$ is singular over the surface of the triangle, the integral of Green's function can be evaluated analytically by employing polar coordinates, see Figure \ref{triangle}. \\
The triangle was splitted into three triangles ($T_1$, $T_2$ and $T_3$. Triangles $T_2$ and $T_3$ are equivalent) considering its centroid as a shared corner. The  corners \textcircled{1}, \textcircled{2}, \textcircled{3} and the centroid of the triangle are described by the coordinates $(0,0)$, $(l,0)$, $(0,l)$ and $(\frac{l}{3},\frac{l}{3})$, respectively. In triangle $T_1$, it is easy to see that $\cos(\alpha)=\frac{3}{\sqrt{10}}$ and $\theta_1 = \pi-2\alpha$. Besides, $\rho_1 = \frac{\sqrt{2}}{6}l$ and $r_1(\theta) = \rho_1\sec(\theta+\alpha-\frac{\pi}{2})$ with $0\le\theta\le\theta_1$ describes the straight line in polar coordinates from \textcircled{2} to \textcircled{3}. Hence, the integral of function $G(\bx,\bxp)$ in triangle $T_1$ is
\begin{eqnarray*}
    \int_{\triangle_{T_1}} G(\bx,\bxp)d\bxp & = & \frac{1}{4\pi}\int_0^{\theta_1}\int_0^{\rho_1}d\rho d\theta = \frac{1}{4\pi}\int_0^{\theta_1}\frac{\sqrt{2}}{6}l\sec\left(\theta+\alpha-\frac{\pi}{2}\right)d\theta \\
    & = & \frac{\sqrt{2}}{12\pi}l\ln{(\sqrt{10}+3)}
\end{eqnarray*}
In the triangle $T_2$, we have $\theta_2 = \frac{3\pi}{4}-\beta=\frac{\pi}{2}+\alpha$ and $\rho_2 = \frac{l}{3}$. Therefore, the straight line in polar coordinates which describes the path from \textcircled{3} to \textcircled{1} is $r_2(\theta)=\rho_2\sec(\theta+\alpha-\frac{5\pi}{4})$, with $\theta_1\le\theta\le\theta_2$. For this reason, the integral of the function $G(\bx,\bxp)$ in the triangle $T_2$ is
\begin{eqnarray*}
    \int_{\triangle_{T_2}} G(\bx,\bxp)d\bxp & = & \frac{1}{4\pi}\int_{\theta_1}^{\theta_1+\theta_2}\int_0^{\rho_2}d\rho d\theta = \frac{1}{4\pi}\int_{\theta_1}^{\theta_1+\theta_2}\frac{1}{3}l\sec\left(\theta+\alpha-\frac{5\pi}{4}\right)d\theta \\
    & = & \frac{1}{12\pi}l\ln|(\sqrt{2}+1)(\sqrt{5}+2)|
\end{eqnarray*}
hence, the integral of function $G(\bx,\bxp)$ in the whole triangle is
\begin{eqnarray*}
    \int_\triangle G(\bx,\bxp)d\bxp & = & \int_{\triangle_{T_1}} G(\bx,\bxp)d\bxp +2\int_{\triangle_{T_2}} G(\bx,\bxp)d\bxp \\
    & = & \frac{l}{12\pi}\ln|(\sqrt{10}+3)^{\sqrt{2}}(\sqrt{5}+2)^2(\sqrt{2}+1)^2|
\end{eqnarray*}
The two dimensional case of integral \eqref{int_green} is well known and is solved in a similar way, see \cite{Brebbia2012BEM}.

\section{Numerical results}

\subsection{A qualitative comparison of the approximation \eqref{Eqn_Approximation}}
\begin{figure}[h]%
    \centering
    \subfloat[\centering Solution of the regularized PDE \eqref{Eqn_RegularisedElliptic}]{{\includegraphics[width=7cm]{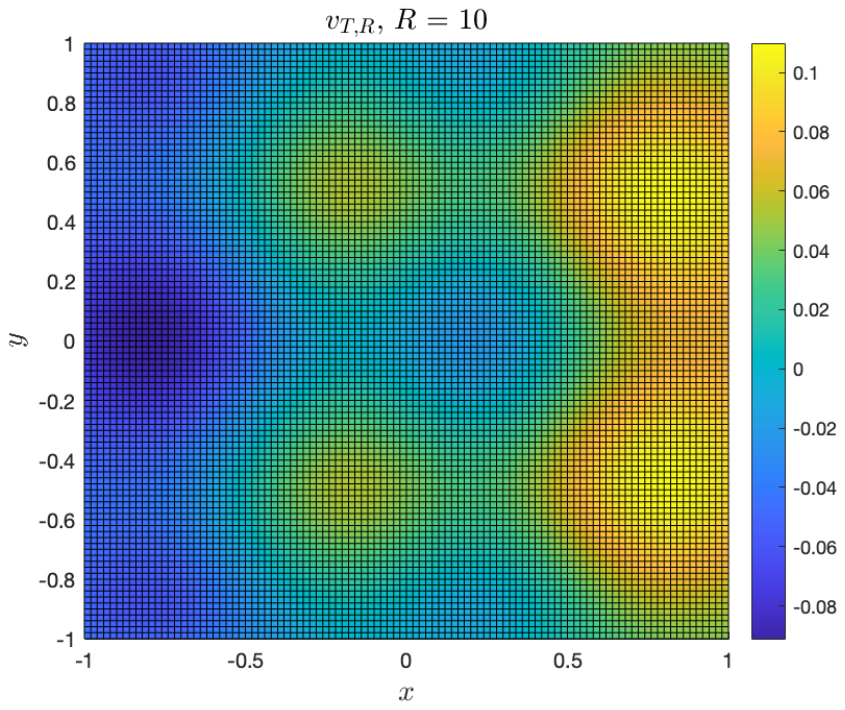} }}%
    \qquad
    \subfloat[\centering Solution of the Volume integral in  \eqref{Eqn_Approximation}]{{\includegraphics[width=7cm]{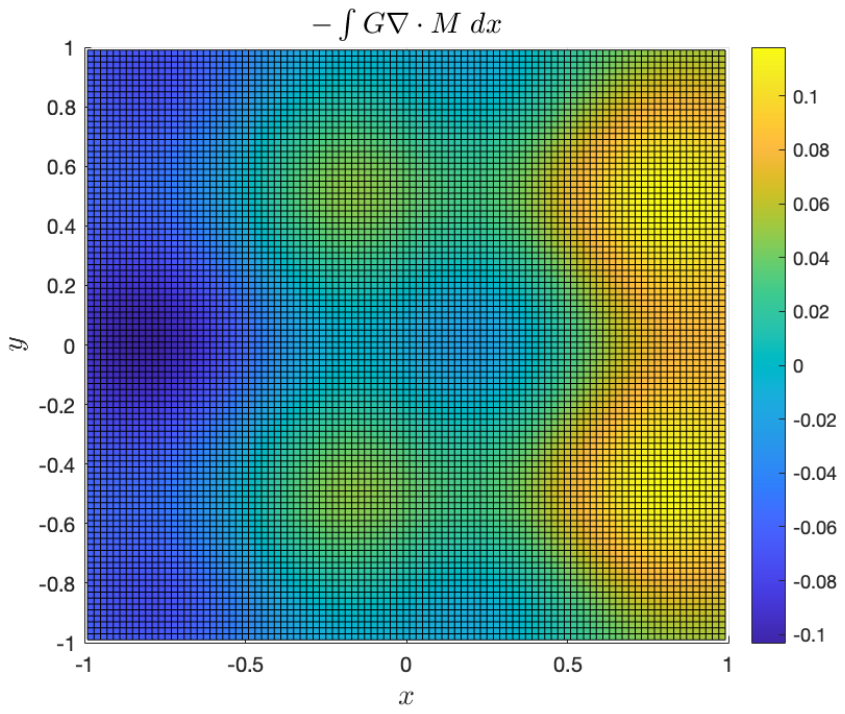} }}%
    \caption{A qualitative comparison of the PDE solution $v_{T,R}$ with the volume integral in \eqref{Eqn_Approximation}}%
    \label{fig:QualitativeComp}%
\end{figure}

Here, our aim is to do a qualitative comparison of $v_{T,R}$ solving \eqref{Eqn_RegularisedElliptic} and the volume integral in \eqref{Eqn_Approximation}. In Figure \ref{fig:QualitativeComp}, both solutions are obtained by using a right hand side 

\begin{equation*}
    F(\bx) = - \nabla \cdot \bM(\bx),  \quad  \bx \in \Omega = [-1,1]^2,  \quad F(\bx) = 0 \quad \bx \in \mathbb{R}^d/\Omega,
\end{equation*}
where 
\begin{align*}
    \bM(\bx)  = \dfrac{1}{2\pi} [\sin(2\pi x_1), -\cos(2\pi x_2)].
\end{align*}

The regularized approach is solved with the parameters $R=10$, and $T = k_T R$, where $k_T \approx 0.18$. Moreover, the number of basis vectors in the Arnoldi algorithm is set to $k = 700$. Both simulations use $N = 50$ points per unit length. In other words, the regularised approach uses $2RN$ points in each direction, whereas the volume integral uses $2N$ discrete points in each direction. Figure \ref{fig:QualitativeComp} demonstrates that both solutions match qualitatively. Note that for the regularised approach, we have included only the restriction of $v_{T,R}$ to $\Omega = [-1,1]^2$ in the figure.

\subsection{Convergence study}
\begin{figure}[h]
    \begin{center}
    \includegraphics[scale=0.5]{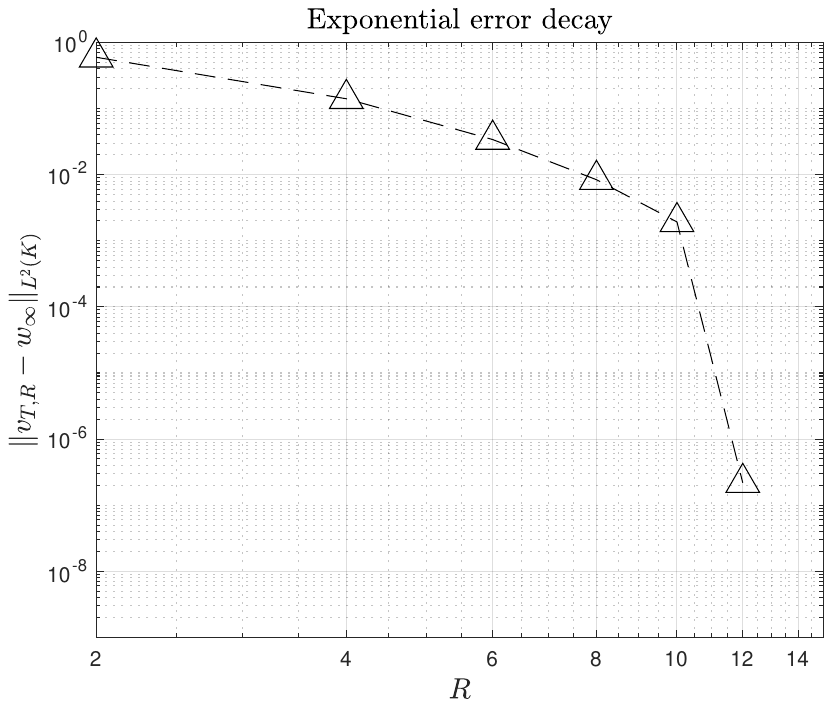}
    \caption{Exponential convergence of $v_{T,R}$ for a periodic problem.}
    \label{fig:PeriodicConvergence}
    \end{center}
\end{figure}

Here, we do a numerical study of convergence of $v_{T,R}$ solving \eqref{Eqn_RegularisedElliptic} to the solution $w_{\infty}$ solving the infinite domain problem \eqref{Eqn_InfDomain_Elliptic}, mainly to demonstrate the convergence rates in Theorems \ref{Thm_Analysis_b} and \ref{Thm_Main_Analysis_c}. For the first test problem, we choose a periodic solution $u(x_1,x_2) = \cos(2 \pi x_1) + \sin(2\pi x_2)$, which implies a right hand side of the form $F(x_1,x_2) = (2\pi)^2 u(x_1,x_2)$. Figure \ref{fig:PeriodicConvergence} depicts the exponential convergence of $v_{T,R}$ to $w_{\infty}$ as $R$ gets larger. Note that in the simulations we have used $T = 0.18 R$, $k = 700$, and $50$ discrete points per unit length (in each direction).

\begin{figure}[h]%
    \centering
    \subfloat[\centering $|\hat{F}(\bome)|$]{{\includegraphics[width=7cm]{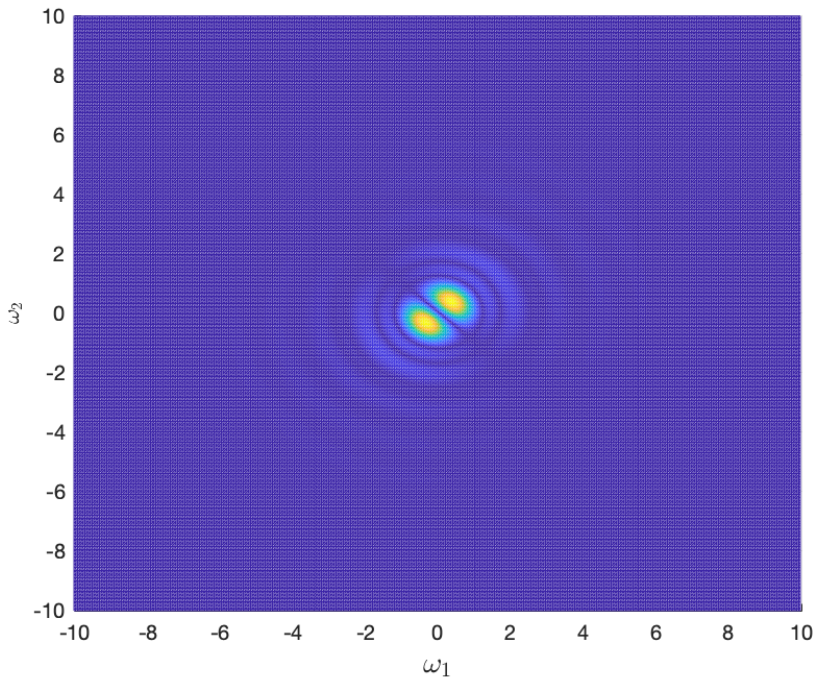} }}%
    \qquad
    \subfloat[\centering $|\hat{F}(\bome)|$ zoomed in the origin ]{{\includegraphics[width=7cm]{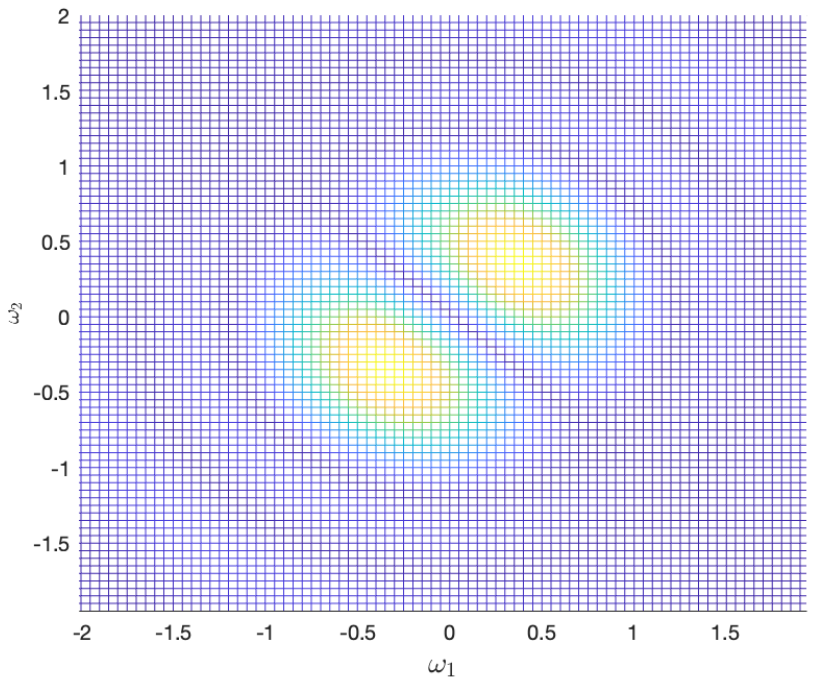} }}%
    \caption{The absolute value of the Fourier transform of F}%
    \label{fig:Fourier_Transform_F}%
\end{figure}

Next, we consider a right hand side $F$ with vanishing frequency components in the origin. We generate the right hand side $F$ as follows:

\begin{align*}
    F(\bx) = \partial_{x_1} \eta(\bx) + \partial_{x_2} \eta(\bx),
\end{align*}
where 
\begin{align*}
    \eta(\bx) = \begin{cases} \exp(\dfrac{1}{| \bx |  - 1}) & |\bx| < 1 \\
    0 & |\bx| \geq 1. 
    \end{cases}
\end{align*}
The absolute of Fourier transform of $F$ is depicted in Figure \ref{fig:Fourier_Transform_F}. From this figure, we visually infer that $\hat{F}(\bome) = 0, \forall \quad |\bome| \leq \omega_0 = \frac{\sqrt{2}}{8}$, which will be needed to choose an optimal value for the parameter $T \approx R/\omega_0$.

\begin{figure}[h]
    \begin{center}
    \includegraphics[scale=0.5]{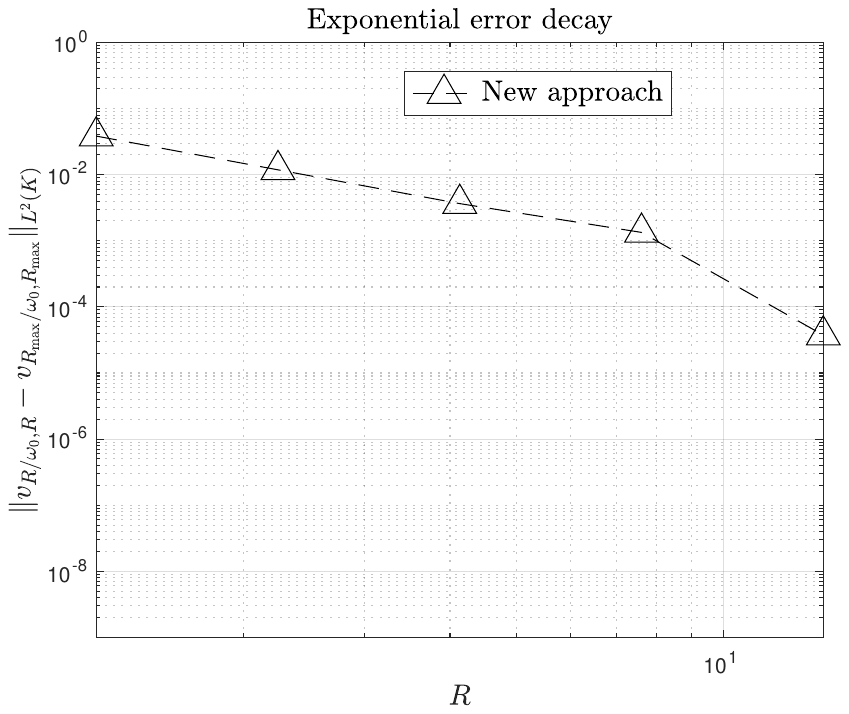}
    \caption{Exponential convergence of $v_{T,R}$ for a problem with vanishing frequency components in the origin ($\omega_0 = \sqrt{2}/8$).}
    \label{fig:HighFreqConvergence}
    \end{center}
\end{figure}

Since for this problem, we do not know the exact solution, we use the solution $v_{T,R}$ of the problem \eqref{Eqn_RegularisedElliptic} with $R= R_{\max}  = 15$, and $T = R_{\max}/\omega_0$ as a reference, and study the convergence of $v_{T,R}$ to the reference solution as the value of the parameter $R$ increases. Figure \ref{fig:HighFreqConvergence} illustrates the exponential convergence (expected from Theorem \ref{Thm_Main_Analysis_c}) with increasing value of $R$. Note that in these simulations, we used $10$ discrete points per unit length in every direction, $k=\min\{700,N\}$, where $N$ is the number of degrees of freedom in one direction, and the value of the parameter $T$ was chosen near optimally as $T = R/\omega_0$, where $\omega_0 = \frac{\sqrt{2}}{8}$.   

\clearpage

\bibliographystyle{plain}
\bibliography{main
}

\end{document}